\documentclass[12pt,a4paper,leqno]{article}

		\title{A spectral sequence for Iwasawa adjoints}
		\author{Uwe Jannsen}
\usepackage{amssymb}
\usepackage{epsf,pstricks,epsfig}
\input xy
\xyoption{all}

\usepackage{epic,eepic}
\usepackage{graphicx}
\newcommand{\ZZ}{\mathbb Z}

\def\noi{\noindent}
\def\bs{\bigskip}
\def\ms{\medskip}

\begin{document}


\maketitle

\section{Introduction}
This paper is a slightly edited and corrected
version of a long time unpublished but several times quoted preprint from 1994.
The aim of this paper was and is to give a purely algebraic tool for treating so-called (generalized) Iwasawa
adjoints of some naturally occurring Iwasawa modules for $p$-adic Lie group extensions,
by relating them to certain continuous Galois cohomology groups via a spectral sequence.

Let $k$ be a number field, fix a prime $p$, and let
$k_\infty$ be some Galois extension of $k$ such that $\mathcal G =
\mbox{Gal }(k_\infty/k)$ is a $p$-adic Lie-group (e.g., $\mathcal
G\cong \mathbb Z_p^r$ for some $r\geqslant 1$). Let $S$ be a
finite set of primes containing all primes above $p$ and $\infty$,
and all primes ramified in $k_\infty/k$, and let $k_S$ be the
maximal $S$-ramified extension of $k$; by assumption, $k_\infty
\subseteq k_S$. Let $G_S= \mbox{Gal }(k_S/k)$ and $G_{\infty,
S}=\mbox{Gal }(k_S/k_\infty)$.
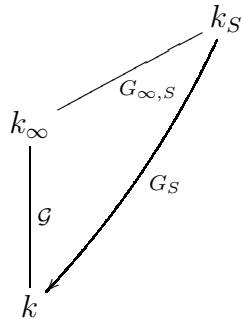
$$
\xymatrix{ & &  k_S\ar@{-}[lld]^{G_{\infty, S}}\ar@/^/[llddd]^{G_S} \\
k_{\infty}\ar@{-}[dd]^{\mathcal G}\\
\\
k}
$$

Let $A$ be a discrete $G_S$-module which is isomorphic to
$(\mathbb Q_p/\mathbb Z_p)^r$ for some $r\geqslant1$ as an abelian
group (e.g., $A=\mathbb Q_p/\mathbb Z_p$ with trivial action, or
$A = E[p^\infty]$, the group of $p$-power torsion points of an
elliptic curve $E/k$ with good reduction outside $S$). We are
\textbf{not} assuming that $G_{\infty, S}$ acts trivially.

\ms\noi Let $\Lambda=\mathbb Z_p[[\mathcal G]]$ be the completed
group ring. For a finitely generated $\Lambda$-module $M$ we put
$$
E^i(M)=\mbox{Ext}^i_\Lambda(M,\Lambda)\; .
$$
Hence $E^0 (M)= \mbox{Hom}_\Lambda(M,\,\Lambda)= :M^+$ is just the
$\Lambda$-dual of $M$. This has a natural structure of a
$\Lambda$-module, by letting $\sigma\in \mathcal G$ act via
$$
\sigma f(m)\; = \; \sigma f(\sigma^{-1})m)
$$
for $f\in M^+$, $m\in M$. It is known that $\Lambda$ is a
noetherian ring (here we use that $\mathcal G$ is a $p$-adic Lie
group), by results of Lazard \cite{La}. Hence $M^+$ is a finitely
generated $\Lambda$-module again (choose a projection
$\Lambda^r\twoheadrightarrow M$; then we have an injection
$M^+\hookrightarrow(\Lambda^r)^+=\Lambda^r$). By standard
homological algebra, the $E^i(M)$ are finitely generated
$\Lambda$-modules for all $i\geqslant 0$ which we call the
(generalized) Iwasawa adjoints of $M$. They can also be seen
as some kind of homotopy invariants of $M$, see \cite{Ja3}, and
also \cite{NSW} V \S 4 and 5.

\ms\noi\textbf{Examples}
(a) If $\mathcal G=\mathbb Z_p$, then $\Lambda=\mathbb
Z_p[[\mathcal G]]\cong\mathbb Z_p[[X]]$ is the classical Iwasawa
algebra, and, for a $\Lambda$-torsion module $M$, it is known (see \cite{P-R} I.2.2, or \cite{Bi} 1.2 and r\'emarque,
or \cite{NSW} Prop. (5.5.6)))
that $E^1(M)$ is isomorphic to the classical Iwasawa adjoint, which was defined by Iwasawa (\cite{Iw} 1.3) as
$$
\mbox{ad}(M)=\lim\limits_{\stackrel {\displaystyle\leftarrow} n}(M/\alpha_n M)^\vee
$$
where $(\alpha_n)_{n\in\mathbb N}$ is any sequence of elements in
$\Lambda$ such that $\lim\limits_{n\rightarrow\infty} \alpha_n=0$ and
$(\alpha_n)$ is prime to the support of $M$ for every $n\geqslant
1$, and where
$$
N^\vee=\mbox{Hom}\,(N,\, \mathbb Q_p/\mathbb Z_p)
$$
is the Pontrjagin dual of a discrete or compact $\mathbb Z_p$-module $N$. For
any finitely generated $\Lambda$-module $M$, $E^1(M)$ is
quasi-isomorphic to $\mbox{Tor}_\Lambda(M)^\sim$, where
$\mbox{Tor}_\Lambda(M)$ is the $\Lambda$-torsion submodule of $M$,
and $M^\sim$ is the ``Iwasawa twist'' of a $\Lambda$-module $M$:
the action of $\gamma\in\mathcal G$ is changed to the action of
$\gamma^{-1}$.

(b) If $\mathcal G=\mathbb Z_p^r$, $r\geqslant 1$, then the
$E^i(M)$ are the standard groups considered in local duality. By
duality for the ring $\mathbb Z_p[[\mathcal G]]=\mathbb
Z_p[[x_1,\ldots ,x_r]]$, they can be computed in terms of local
cohomology groups (with support) or by a suitable Koszul complex.
More precisely, $E^i(M) \cong H^{r+1-i}_{\mathfrak{m}}(M)^\vee$,
as recalled in \cite{Bi} 1.

\ms\noindent The main result of this note is the following observation.

\bs\noi\textbf{Theorem 1} \emph{There is a spectral sequence of
finitely generated $\Lambda$-modules}
$$
E_2^{p,q} \!\!=\! E^p(H^q(G_{\infty,S},
A)^\vee)\!\Rightarrow\!\lim\limits_{\stackrel {\displaystyle\leftarrow} {k',m}} H^{p+q}(G_S(k'),
A[p^m])\!\!=\!\lim\limits_{\stackrel {\displaystyle\leftarrow} {k'}} H^{p+q}(G_S(k'),\,T_pA) \,
.
$$
\emph{Here the inverse limits runs through the finite extensions $k'/k$ contained in $k_\infty$,
and the natural numbers $m$,
via the corectrictions and the natural maps
$$
H^n(G_S(k'),\, A[p^{m+1}])\rightarrow H^n(G_S(k'),\, A[p^m])\,,
$$
respectively. The groups
$$
H^{p+q}(G_S(k'),\,T_pA) = \lim\limits_{\stackrel {\displaystyle\leftarrow} {m}} H^{p+q}(G_S(k'),\,A[p^m])
$$
are the continuous cohomology groups of the Tate module} $T_pA =
\lim\limits_{\stackrel {\displaystyle\leftarrow} m} A[p^m]$\,.

\section{Some consequences}

Before we give the proof of a slightly more general result
(cf. Theorem 11 below), we discuss what this spectral sequence
gives in more down-to-earth terms. First of all, we always have
the 5-low-terms exact sequence
$$
0  \longrightarrow  E^1(H^0(G_{\infty, S},A)^\vee)
\stackrel{\scriptstyle \mbox{inf}^1}{\longrightarrow}
\lim\limits_{\stackrel {\displaystyle\leftarrow} {k'}} H^1(G_S(k'), T_p A)
\phantom {0}\longrightarrow  (H^1(G_{\infty,\, S},\, A)^\vee)^+
$$
$$
\longrightarrow  E^2(H^0(G_{\infty, S}, A)^\vee)
\stackrel{\scriptstyle \mbox{inf}^2}{\longrightarrow}
\lim\limits_{\stackrel {\displaystyle\leftarrow} {k'}} H^2(G_S(k'), T_p A) .
$$
To say more, we make the following assumption.

\bs\noi \textbf{A.1} Assume that $p>2$ or that $k_\infty$ is totally imaginary.

\bs\noi
It is well-known that this implies
$$
\leqno{\textbf{A.2}}\;\;
H^r(G_{\infty,\, S},\, A)=0=\lim\limits_{\stackrel {\displaystyle\leftarrow} {k'}} H^r(G_S(k'),T_p A) \; \mbox{  for all  }r>2\,.
$$

\smallskip\noi\textbf{Corollary 2} \emph{Assume in addition that $H^2(G_{\infty,\,
S}, \, A)=0$. (This is the so-called ``weak Leopoldt conjecture''
for $A$. It is stated classically for $A=\mathbb Q_p/\mathbb Z_p$
with trivial action, and there are precise conjectures when this is expected to hold for
modules $A$ coming from algebraic geometry, cf. \cite{Ja2}.)
Then the cokernel of $\inf^2$ is
$$
\ker\,(E^1(H^1(G_{\infty,\, S},\, A)^\vee)\twoheadrightarrow
E^3(H^0(G_{\infty,\, S},\, A)^\vee))\, ,
$$
and there are isomorphisms
$$
E^i(H^1(G_{\infty,\, S},\,
A)^\vee) \stackrel{\sim}{\longrightarrow}
E^{i+2}(H^0(G_{\infty,\, S},\, A)^\vee)
$$
for $i\geqslant 2$.}

\bs\noi\textbf{Proof} This comes from A.2 and the following
picture of the spectral sequence
$$
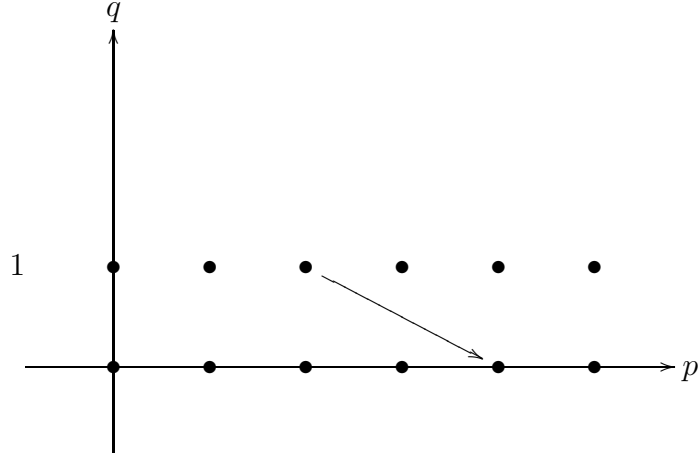
\xymatrix{ & q\\
\\
\\
1 & \bullet & \bullet & \bullet\ar[drr] & \bullet & \bullet & \bullet\\
\ar[rrrrrrr]& \bullet & \bullet & \bullet & \bullet & \bullet & \bullet & p\\
& \ar[uuuuu] }
$$

\bs\noi\textbf{Corollary 3} \emph{Assume that $H^0(G_{\infty,\,
S},\, A)=0$. Then}\\
\noi (a)
$$
\lim\limits_{\stackrel {\displaystyle\leftarrow} {k'}} H^1(G_S(k'),\, T_pA)\stackrel{\sim}{\longrightarrow}
H^1(G_{\infty,\, S},\, A)^\vee)^+\;.
$$
\\
\noi (b) \emph{ There is an exact sequence}
$$
0\longrightarrow E^1(H^1(G_{\infty,\, S},\, A)^\vee)\longrightarrow
\lim\limits_{\stackrel {\displaystyle\leftarrow} {k'}} H^2(G_S(k'),\, T_pA)
$$
$$
\longrightarrow ( H^2(G_{\infty,\, S},\, A)^\vee)^+\longrightarrow E^2
(H^1(G_{\infty,\, S},\, A)^\vee)\longrightarrow 0\;.
$$
\noi (c) \emph{There are isomorphisms}
$$
E^i (H^2(G_{\infty,\, S},\,
A)^\vee)\stackrel{\sim}{\longrightarrow} E^{i+2}( H^1(G_{\infty,\,
S},\, A)^\vee)
$$
\emph{for $i\geqslant 1$.}

\bs\noi\textbf{Proof} In this case, the spectral sequence looks like
$$
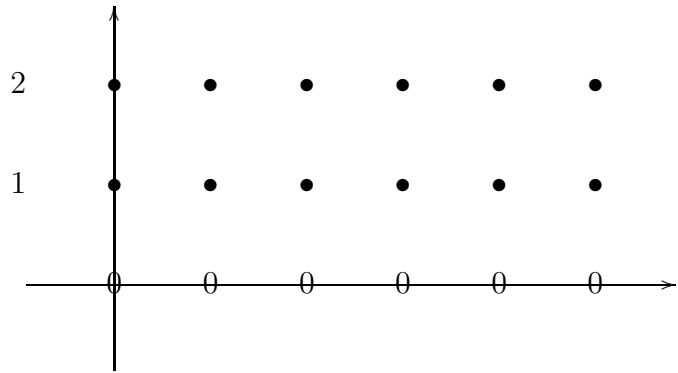
\xymatrix{
&\\
2 & \bullet & \bullet & \bullet & \bullet & \bullet & \bullet\\
1 & \bullet & \bullet & \bullet & \bullet & \bullet & \bullet\\
\ar[rrrrrrr]& 0 & 0 & 0 & 0 & 0 & 0 & \\
& \ar[uuuu] }
$$

\bs\noi\textbf{Corollary 4} \emph{Assume that $\mathcal G$ is a
$p$-adic Lie group of dimension $1$ (equivalently: an open
subgroup is $\cong\mathbb Z_p$). Then $E^i(-)=0$ for $i\geqslant
3$. Let}
$$
B= \mbox{im }(\inf{^2}:E^2(H^0(G_{\infty,\, S},\,
A)^\vee)\rightarrow\lim\limits_{\stackrel {\displaystyle\leftarrow} {k'}} H^2(G_S(k'),\, T_p A))
$$
\emph{Then $B$ is finite, and there is an exact sequence}
$$
0 \rightarrow  E^1(H^1(G_{\infty,S},\, A)^\vee)\rightarrow
\lim\limits_{\stackrel {\displaystyle\leftarrow} {k'}} H^2(G_S(k'), T_p A)/B
 \rightarrow (H^2(G_{\infty,
S}, A)^\vee)^+
$$
$$
\rightarrow E^2(H^1(G_{\infty,\, S},\, A)^\vee)\rightarrow 0\, ,
$$
\emph{and}
$$
E^1(H^2(G_{\infty,\, S},\, A)^\vee)=0= E^2(H^2(G_{\infty,\, S},\,
A)^\vee)\, ,
$$
\emph{i.e., $(H^2(G_{\infty,\, S},\, A)^\vee$ is a projective
$\Lambda$-module.}

\bs\noi\textbf{Proof} Quite generally, for a $p$-adic Lie group
$\mathcal G$ of dimension $n$ one has $vcd_p(\mathcal G)=n$ for
the virtual cohomological $p$-dimension of $\mathcal G$, and hence
$E^i(-)=0$ for $i>n+1$, cf. \cite{Ja3} Cor. 2.4. The finiteness of $E^2(M)$,
for a $\Lambda$-module $M$ which is finitely generated over $\mathbb Z_p$
(like our module $H^0(G_{\infty,\, S},\,A)^\vee$) follows from Lemma 5 below.
In fact, the exact sequence $0 \rightarrow M_{tor} \rightarrow M \rightarrow \tilde{M} \rightarrow 0$,
in which $M_{tor}$ is the torsion submodule of $M$, induces a long exact sequence
$$
\ldots \rightarrow E^i(\tilde{M}) \rightarrow E^i(M) \rightarrow E^i(M_{tor}) \rightarrow \ldots
$$
in which we have $E^i(M_{tor}) = 0$ for $i\neq n+1$ and finiteness of $E^{n+1}(M_{tor})$
by Lemma 5 (b), and $E^i(\tilde{M})=0$ for $i \neq n$ for the torsion-free module $\tilde{M}$
by Lemma 5 (a). In our case we have $n=1$ and therefore the finiteness of $E^2(M)$.
The remaining claims follow from the following shape of the spectral sequence:
$$
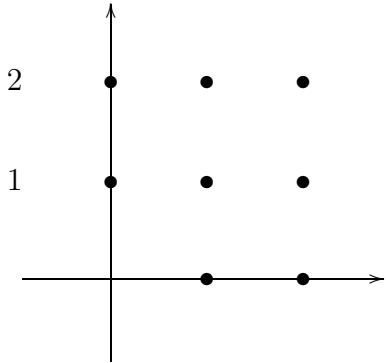
\xymatrix{
& & \\
2 & \bullet & \bullet & \bullet \\
1 & \bullet & \bullet & \bullet \\
\ar[rrrr]&  & \bullet & \bullet &\\
& \ar[uuuu] }
$$

\bs\noi\textbf{Lemma 5} \emph{Assume that $\mathcal G$ is a
$p$-adic Lie group of dimension $n$ (this holds, e.g., if
$\mathcal G$ contains an open subgroup $\cong\mathbb Z_p^n$),
and let $M$ be a $\Lambda$-module which is finitely generated
as a $\mathbb Z_p$-module. Then the follwing holds.}\\
(a) $E^i (M)= 0$ for $i\neq n,\, n+1$.\\
(b) \emph{If $M$ is torsion-free, then}
$$
E^i(M)=\left\{\begin{array}{ccc}0
& \mbox{for} & i\neq n\\
\mbox{Hom } \,(D,\, M^\vee),& \mbox{for} & i =
n\, ,
\end{array}\right.
$$
\emph{where $D$ is the dualising module for $\mathcal G$ (which is a
divisible cofinitely generated $\mathbb Z_p$-module, e.g., $D=\mathbb Q_p/\mathbb Z_p$ if $\mathcal G=\mathbb Z_p^n$).
In particular, $E^n(M)$ is a torsion-free finitely generated $\mathbb Z_p$-module.}\\
(c) \emph{If $M$ is finite, then}
$$
E^i(M)=\left\{\begin{array}{ccc}0
& \mbox{for} & i\neq n+1\\
\mbox{Hom } \,(M,\, D)^\vee & \mbox{for} &
i = n+1
\end{array}\right.
$$
\emph{In particular, $E^{n+1}(M)$ is a finite $\mathbb Z_p$-module.}

\bs\noi\textbf{Proof} See \cite{Ja3}, Cor. 2.6. For (a) note the isomorphism
$$
\mbox{Hom}(D\otimes_{\mathbb Z_p}M,\mathbb Q_p/\mathbb Z_p) \cong \mbox{Hom}(D,\mbox{Hom}(M,\mathbb Q_p/\mathbb Z_p))
= \mbox{Hom}(D,M^\vee)\;.$$

\bs\noi\textbf{Corollary 6} \emph{Let $\mathcal G$ be a $p$-adic
Lie group of dimension 2 (e.g., $\mathcal G$ contains an open
subgroup $\cong \mathbb Z^2_p$). If $G_{\infty,\, S}$ acts
trivially on $A$, then there are exact sequences}
{\normalsize
$$
0\rightarrow\lim\limits_{\stackrel {\displaystyle\leftarrow} {k'}} H^1(G_S(k'),\, T_p A)\rightarrow(H^1(G_{\infty,\, S},\,
A)^\vee)^+
\rightarrow
T_pA\stackrel{\inf^2}{\longrightarrow}\lim\limits_{\stackrel {\displaystyle\leftarrow} {k'}} H^2(G_S(k'),\, T_p A)
$$
}
\emph{and}
$$
\begin{array}{rl}
 0\rightarrow E^1(H^1(G_{\infty,\, S},\,
A)^\vee)\rightarrow\lim\limits_{\stackrel {\displaystyle\leftarrow} {k'}} H^2(G_S(k'),\, T_p A)
/\mbox{im} \inf^{2}
\end{array}
$$
$$
\rightarrow(H^1(G_{\infty,\, S},\, A)^\vee)^+\\
\rightarrow E^2(H^1(G_{\infty,\, S},\, A)^\vee)\rightarrow 0\, ,
$$
\emph{an isomorphism}
$$
E^1(H^2(G_{\infty,\, S},\,
A)^\vee)\stackrel{\sim}{\longrightarrow} E^3 (H^1(G_{\infty,\,
S},\, A)^\vee),\,
$$
\emph{and one has}
$$
E^2(H^2(G_{\infty,\, S},\, A)^\vee)\;=\; 0\; =\; E^3
(H^2(G_{\infty,\, S},\, A)^\vee)\, .
$$

\bs\noi\textbf{Proof} The spectral sequence looks like
$$
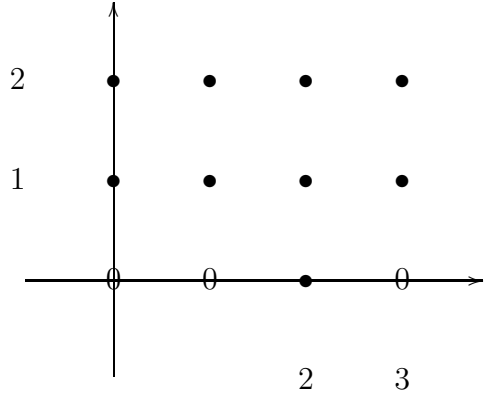
\xymatrix{
&\\
2 & \bullet & \bullet & \bullet & \bullet \\
1 & \bullet & \bullet & \bullet & \bullet \\
\ar[rrrrr]& 0 & 0 & \bullet & 0 & \\
& \ar[uuuu] & & 2 & 3  }
$$

\bs\noi\textbf{Corollary 7} \emph{Let $\mathcal G$ be a $p$-adic
Lie group of dimension 2 (So $E^i(-)=0$ for $i\geqslant 4$). If
$H^0(G_{\infty,\, S},\, A)$ is finite, then}
$$
\lim\limits_{\stackrel {\displaystyle\leftarrow} {k'}} H^1(G_S(k'),\, T_p A)\cong (H^1(G_{\infty,\, S},\,
A)^\vee)^+\; .
$$
\emph{If}
$$
d_2^{1,1}: E^1(H^1(G_{\infty,\, S},\, A)^\vee)\; \longrightarrow\;
E^3(H^0(G_{\infty,\, S},\, A)^\vee)
$$
\emph{is the differential of the spectral sequence in the theorem,
then one has an exact sequence}
$$
0\;\longrightarrow \; \ker d_2^{1,1} \; \longrightarrow
\lim\limits_{\stackrel {\displaystyle\leftarrow} {k'}} H^2(G_S(k'),\, T_p A)\; \longrightarrow
$$
$$
\longrightarrow \; \ker (d^{0,2}_2\!\!:(H^2(G_{\infty,\, S},\,
A)^\vee)^+\twoheadrightarrow E^2(H^1(G_{\infty,\, S},\, A)^\vee))
\longrightarrow\; \mbox{coker } d^{1,1}_2\longrightarrow 0\,,
$$
\emph{an isomorphism}
$$
E^1(H^2(G_{\infty,\, S},\,
A)^\vee)\stackrel{\sim}{\longrightarrow}E^3(H^1(G_{\infty,\, S},\,
A)^\vee)\, ,
$$
\emph{and the vanishing}
$$
E^2(H^2(G_{\infty,\, S},\, A)^\vee)\;=\;0\;=\;E^3(H^2(G_{\infty,\,
S},\, A)^\vee)\; .
$$

\bs\noi\textbf{Proof} The spectral sequence looks like
$$
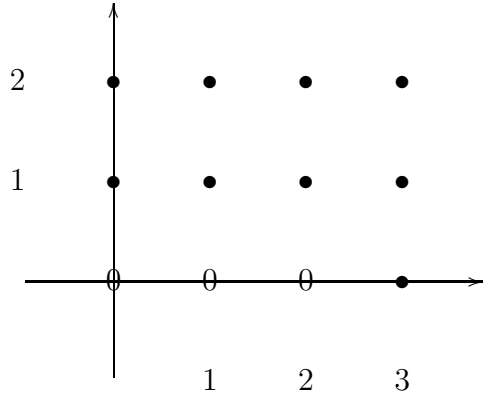
\xymatrix{
&\\
2 & \bullet & \bullet & \bullet & \bullet \\
1 & \bullet & \bullet & \bullet & \bullet \\
\ar[rrrrr]& 0 & 0 & 0 & \bullet & \\
& \ar[uuuu] & 1 & 2 & 3  }
$$

\bs\noi\textbf{Remark} In the situation of Corollary 5, one has an
exact sequence up to \emph{finite modules}:
$$
0\rightarrow E^1(H^1(G_{\infty,\, S},\,
A)^\vee)\rightarrow\lim\limits_{\stackrel {\displaystyle\leftarrow} {k'}} H^2(G_S(k'),\, T_p A)
$$
$$
\rightarrow(H^2(G_{\infty,\, S},\, A)^\vee)^+\rightarrow
E^2(H^1(G_{\infty,\, S},\, A)^\vee)\rightarrow 0\, .
$$

\bs\noi\textbf{Corollary 8} \emph{Let $\mathcal G$ be a $p$-adic
Lie group of dimension $>2$. Then}
$$
(H^1(G_{\infty,\, S},\, A)^\vee)^+\; \cong\;
\lim\limits_{\stackrel {\displaystyle\leftarrow} {k'}} H^1(G_S(k'),\, T_p A)
$$

\bs\noi\textbf{Proof} The first three columns of the spectral
sequence look like
$$
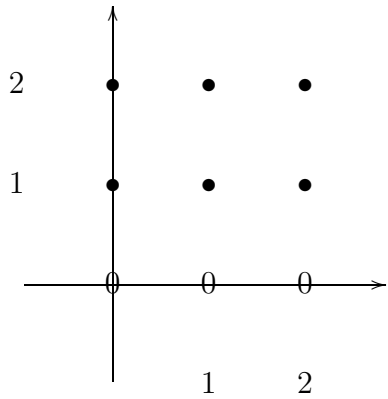
\xymatrix{
&\\
2 & \bullet & \bullet & \bullet  \\
1 & \bullet & \bullet & \bullet  \\
\ar[rrrr]& 0 & 0 & 0 & \\
& \ar[uuuu]  & 1 & 2  }
$$

\section{Proof of the Main Theorem}

\bs\noi We will now prove Theorem 1, by proving a somewhat more
general result. For any profinite group $G$, let $\Lambda(G) =
\mathbb Z_p[[G]]$ be the completed group ring over $\mathbb Z_p$,
and let $M_G = M_{G,p}$ be the category of discrete (left)
$\Lambda(G)$-modules. These are the discrete $G$-modules $A$ which
are $p$-primary torsion abelian groups. For such a module $A$, its
Pontrjagin dual $A^\vee = Hom(A,\mathbb Q_p/\mathbb Z_p)$ is a
compact $\Lambda(G)$-module. In fact, Pontrjagin duality gives an
anti-equivalence between $M_G$ and the category $C_G = C_{G,p}$ of
compact (right) $\Lambda(G)$-modules.

\ms\noi Let $M_G^{\mathbb N}$ be the category of inverse systems
$$(A_n) : ~~\ldots \rightarrow A_3 \rightarrow A_2 \rightarrow A_1
$$
in $M_G$ as in \cite{Ja1}. Denote by $H^i_{cont}(G,(A_n))$ the
continuous cohomology of such a system and recall that one has an
exact sequence for each $i$
$$
0~ \rightarrow~ R^1\mathop{\mbox{lim}}\limits_{\stackrel {\displaystyle\leftarrow} {n}} H^{i-1}(G, A_n)~ \rightarrow~
H^i_{cont}(G,(A_n))~ \rightarrow \lim\limits_{\stackrel {\displaystyle\leftarrow} {n}} H^i(G, A_n)~ \rightarrow~
0 \, ,
$$
in which the first derivative $R^1\mathop{\mbox{lim}}\limits_{\stackrel {\displaystyle\leftarrow} {n}}$ of the inverse limit,
also noted as ${\mathop{\mbox{lim}}\limits_{\stackrel {\displaystyle\leftarrow} {n}}}^1$,
vanishes if the groups
$H^{i-1}(G, A_n)$ are finite for all $n$ (cf. loc. cit.).

\bs\noi\textbf{Definition 9} \emph{For a closed subgroup $H \leq
G$ and a discrete $G$-module $A$ in $M_G$ define the relative
cohomology $H^m(G,H; A)$ as the value at $A$ of the $m$-th derived
functor of the left exact functor (with $Ab$ being the category of
abelian groups)
$$
\begin{array}{rcl}
H^0(G,\, H; \, -):M_G &
\rightarrow & Ab\\
A & \mapsto & \lim\limits_{\stackrel {\displaystyle\leftarrow} {U}} H^0(U,\, A)\, ,
\end{array}
$$
where $U$ runs through all open subgroups $U\subset G$ containing
$H$, and the transition maps are the corestriction maps. For an
inverse system $(A_n)$ of modules in $M_G$ define the continuous
relative cohomology $H^m_{cont}(G,H; (A_n))$ as the value at
$(A_n)$ of the $m$-th right derivative of the functor
$$
\begin{array}{rcl}
H^0_{\mbox{\scriptsize cont}}(G,\, H; \, -):M_G^{\mathbb N} &
\rightarrow & Ab\\
(A_n) & \mapsto &
\lim\limits_{\stackrel {\displaystyle\leftarrow} {n}}\lim\limits_{\stackrel {\displaystyle\leftarrow} {U}}H^0(U,\, A_n)
\, ,
\end{array}
$$
where the limit over $U$ is as before, and the limit over $n$ is
induced by the transition maps $A_{n+1} \rightarrow A_n$.}

\bs\noi\textbf{Lemma 10} \emph{If $G/H$ has a countable basis of
neighbourhoods of identity, i.e., if there is a countable family
$U_\nu$ of open subgroup, $H\leq U_\nu\leq G$, with
$\raisebox{-6pt}{$\stackrel{\bigcap}{\nu}$}\, U_\nu=H$, and if, in
addition, $H^i(U, A_n)$ is finite for all these $U$ and all $n$,
then
$$
H^i_{\mbox{\scriptsize cont}}\,(G,\,
H;\,(A_n))=\lim\limits_{\stackrel {\displaystyle\leftarrow} {n}}\lim\limits_{\stackrel {\displaystyle\leftarrow} {U}}H^i(U,\,A_n)\, .
$$}

\bs\noi\textbf{Proof} More generally, without assuming the finiteness of the groups $H^i(U,A)$,
we claim that we have a Grothendieck spectral sequence for the composition of the functors
$(A_n)_n \rightsquigarrow (H^0(U,A_n))_{U,n}$ with the functor
$\lim\limits_{\footnotesize \leftarrow n,U}$
$$
E_2^{p,q}=
R^p\lim\limits_{\footnotesize \leftarrow n,U}\,H^q(U,\, A_n)\Rightarrow H_{\mbox{\scriptsize
cont}}^{p+q}(G,\, H;\, (A_n))\, .
$$
For this we have to show that the first functor sends injective objects to acyclics for
the second functor. But if $(I_n)$ is an injective system, then all $I_n$ are injective and
all morphisms $I_{n+1} \rightarrow I_n$ are split surjections, see \cite{Ja1} (1.1),
so $H^0(U,I_{n+1}) \rightarrow H^0(U,I_n)$ is surjective for any open subgroup $U\subset G$.

On the other hand, if $I$ is an injective $G$-module, then for any pair of open subgroups
$U' \subset U$ the corestriction $cor: H^0(U',I) \rightarrow H^0(U,I)$ is surjective.
In fact, we may assume that $I = Ind_1^G(B)$ is an induced module for a divisible abelian group $B$.
(Any such module is injective, and any discrete $G$-module can be embedded into such a module,
see \cite{Sha} p. 28 and 29., so that any injective is a direct factor of such a module).
Moreover, since the formation of corestrictions is transitive, we may consider an open subgroup
$U'' \subset U'$ which is normal in $G$. Then $Ind_1^G(B)^{U''} \cong Ind_1^{\overline G}(B)$
for $\overline{G}=G/U''$, and it is known that this is a cohomologically trivial $\overline{G}$-module.
Therefore, letting $\overline{U}=U/U''$, we have
$$
Ind_1^G(B)^U = Ind_1^{\overline{G}}(B)^{\overline{U}} = tr_{\overline{U}}Ind_1^{\overline G}(B) = cor_{U''/U} Ind_1^G(B)^{U''}
$$
as claimed.

By assumption, the inverse limit over the open subgroups $U$ containing $H$ can be replaced by a
cofinal set of subgroups $U_m$ with $m\in\mathbb N$ and $U_{m+1} \subset U_m$, and then
the limit over these $U_m$ and over $n$ can be replaced by the `diagonal' limit over the
pairs $(U_n,n)$ for $n\in \mathbb N$. For such an inverse limit it is well-known that
$R^p\lim\limits_{\leftarrow, n} =0$ for $p>1$, and that $R^1\lim\limits_{\leftarrow, n} H^0(U_n,I_n) = 0$,
since the transition maps $H^0(U_{n+1},I_{n+1}) \rightarrow H^0(U_n,I_n)$ are surjective, as shown above
(so the system trivially satisfies the Mittag-Leffler condition).

This shows the existence of the above spectral sequence.
If, in addition, all $H^q(U,\, A_n)$ are finite, then, reasoning as above,
$$
R^1\lim\limits_{\footnotesize \leftarrow n,U}\,H^q(U,\,A_n)= 0
$$
by the Mittag-Leffler property, and we get the claimed isomorphisms.

\bs\noi Now we come to the spectral sequence in theorem 1. Any
module $A$ in $M_G$ gives rise to two inverse systems, viz., the
system $(A[p^n])$, where the transition maps
$A[p^{n+1}]\rightarrow A[p^n]$ are induced by multiplication with
$p$ in $A$, and the system $(A/p^n)$, where the transition maps
are induced by the identity of $A$. For reasons explained later,
denote by $H^m_{\mbox{\scriptsize cont}}(G,\,H;\,R\underline
T_pA)$ the value at $A$ of the $m$-th derived functor of the left
exact functor
$$
F:A\shortmid\hspace{-0,8mm}\rightsquigarrow\lim\limits_{\stackrel {\displaystyle\leftarrow} {n}}\lim\limits_{\stackrel {\displaystyle\leftarrow} {U}} H^0(U,\,
A[p^n])
$$
where $U$ runs through all open subgroups $U\subset G$ containing
$H$, and the transition maps are the corestriction maps and those
coming from $A[p^{n+1}]\rightarrow A[p^n]$, respectively. If $H$
is a normal subgroup, then we may restrict to normal open
subgroups $U \leq G$ containing $H$ in the above inverse limit,
and the limit is a (left) $\Lambda(G/H)$-module in a natural way.

\bs\noi\textbf{Theorem 11} \emph{Let $H$ be a closed subgroup of a
profinite group $G$ such that $G/H$ has a countable basis of
neighbourhoods of identity (see Lemma 10), and let $A$ be a discrete
$\Lambda(G)$-module.\\
\noi (a)
 There are short exact sequences
 $$
 0\rightarrow H^i_{\mbox{\scriptsize cont}}(G, H;
 (A[p^n]))\rightarrow H^i_{\mbox{\scriptsize cont}}(G, H;
 R\underline T_pA)\rightarrow H^{i-1}_{\mbox{\scriptsize cont}} (G, H;
 (A/p^n))\rightarrow 0 \, .
 $$
 If $H$ is a normal subgroup, then these are exact sequences of
 $\Lambda(G/H)$-modules.\\
 \noi (b)
 Let $H'$ be a normal subgroup of $G$, with $H' \subset H$. There is a spectral sequence
 $$
 E_2^{p,q} = H^p_{cont}(G/H', H/H';\, R\underline T_pH^q(H',\, A)) ~\Rightarrow~
 H_{cont}^{p+q} \,(G,\, H;\, R\underline T_pA)\, .
 $$
 If $H$ is a normal subgroup, too, this is a spectral sequence of
 $\Lambda(G/H)$-modules.\\
 \noi (c)
 If $H$ is a normal subgroup of $G$, then for every discrete $\Lambda(G)$-module
 $A$ one has canonical isomorphisms of $\Lambda(G/H)$-modules
 $$
 H^m(G,H; R\underline T_pA)~\cong~Ext^m_{\Lambda(G)}(A^\vee, \Lambda(G/H))
 $$
 for all $m\geq 0$, where $\Lambda(G/H)$ is regarded as a $\Lambda(G)$-module via the
 ring homomorphism $\Lambda(G) \rightarrow \Lambda(G/H)$. More precisely,
 the $\delta$-functor
 $$
 M_G~\rightarrow ~Mod_{\Lambda(G/H)}~~~,~~~A ~\rightsquigarrow~
 (H^m(G,H; R\underline T_pA) ~|~m\geq 0)
 $$
 is canonically isomorphic to the $\delta$-functor
 $$
 M_G~\rightarrow ~Mod_{\Lambda(G/H)}~~~,~~~A ~\rightsquigarrow~
 (Ext^m_{\Lambda(G)}(A^\vee,\Lambda(G/H))~|~ m\geq 0)\, .
 $$
 Here and in the following, the $Ext$-groups $Ext_{\Lambda(G)}(-,-)$ are
 taken in the category $C_G$ of compact $\Lambda(G)$-modules. We
 note that these $Ext$-groups are $\Lambda(G)$-modules, but not
 necessarily compact.\\
 \noi (d) In particular, let $H$ be a normal subgroup of $G$,
 and let $\mathcal G = G/H$. If $A$ is a discrete $\Lambda(G)$-module,
 then one has a spectral sequence of $\Lambda(\mathcal G)$-modules
 $$
 E_2^{p,q}\! = Ext^p_{\Lambda(\mathcal G)}(H^q(H,\, A)^\vee,
 \Lambda(\mathcal G))\!\Rightarrow \!H_{cont}^{p+q} (G, H;
 R\underline T_pA) \! = Ext^{p+q}_{\Lambda(G)}(A^\vee,\Lambda(\mathcal G)) .
 $$
 }

\bs\noi Before we give the proof of Theorem 11, we note that it
implies Theorem 1. In fact, we apply Theorem 11 to $G = G_S$ and
$H = G_{\infty,S}$. If $A$ is a $G_S$-module of cofinite type as
in Theorem 1, then $A/p^n = 0$ and $A[p^n]$ is finite, for all
$n$. Moreover, $H^i(U,B)$ is known to be finite for all open
subgroups $U \leq G_S$ and all finite $U$-modules $B$. By (a) and
Lemma 10 we deduce
$$
H_{cont}^m \,(G_S,\, G_{\infty,S};\, RT_pA) = \lim\limits_{\stackrel {\displaystyle\leftarrow} {n,U}}  H^m(U,
A[p^n]) = \lim\limits_{\stackrel {\displaystyle\leftarrow} {n,k'}}H^m(G_S(k'), A[p^n])\, ,
$$
where $k'$ runs through all finite subextensions of $k_\infty/k$.
Moreover, one has canonical isomorphisms
$$
\lim\limits_{\stackrel {\displaystyle\leftarrow} {n}} H^m(U, A[p^n])~\cong~H^m(U,T_pA)
$$
where the latter group is continuous cochain group cohomology, cf.
[Ja 1]. By applying Theorem 11 (d) we thus get the desired
spectral sequence. Finally, $H^m(H,A)^\vee$ is a finitely
generated $\Lambda(\mathcal G)$-module for all $m\geq 0$, so that
the initial terms of the spectral sequence are finitely generated
$\Lambda(\mathcal G)$-modules as well, and so are the limit terms.
In fact, let $N$ be the kernel of the homomorphism $G_S
\rightarrow Aut(A)$ given by the action of $G_S$ on $A$, and let
$H' = H\cap N$. Then $G/H'$ is a $p$-adic analytic Lie group,
since $G/H$ and $G/N$ are. It is well-known that $H^m(H',\mathbb
Q_p/\mathbb Z_p)$ is a cofinitely generated discrete
$\Lambda(G/H')$-module for all $m\geq 0$; hence the same is true
for $H^m(H',A) \cong H^m(H',\mathbb Q_p/\mathbb Z_p)\otimes T_pA$.
The claim then follows from the Hochschild-Serre spectral sequence
$H^p(H/H',H^q(H',A)) \Rightarrow H^{p+q}(H,A)$.

\ms\noi\textbf{Proof of Theorem 11} (a): We can write $F$ as the
composition of the two left exact functors
$$
\begin{array}{rcl}
\underline T_p:M_G & \rightarrow & M_G^{\mathbb N}\\
A & \shortmid\hspace{-0,8mm}\rightsquigarrow & (A[p^n])
\end{array}
$$
and
$$
\begin{array}{rcl}
H^0_{\mbox{\scriptsize cont}}(G,\, H; \, -):M_G^{\mathbb N} &
\rightarrow & Ab\\
(A_n) & \rightsquigarrow &
\lim\limits_{\stackrel {\displaystyle\leftarrow} {n}}\lim\limits_{\stackrel {\displaystyle\leftarrow} {U}} H^0(U,\, A_n) \, ,
\end{array}
$$
where the limit over $U$ runs through all open (normal) subgroups of $G$
containing $H$, with the corestrictions as transition maps. With the arguments in the proof of lemma 10, we can deduce that
$\underline T_p$ maps injectives to $H^0_{\mbox{\scriptsize cont}}(G,\, H;\, -)$-acyclics.
In fact, we may assume injective $G$-modules given as induced modules $I = Ind_1^G(B)$ with a divisible
abelian group $B$. Then each module $I[p^n] = Ind_1^G(B[n])$ is induced, hence acyclic for the functor $H^0(U,-)$,
and since $I$ is divisible, we have exact sequences
$$
0 \rightarrow I[p] \longrightarrow I[p^{r+1}] \longrightarrow I[p^r] \rightarrow 0\,.
$$
Therefore the transition maps $H^0(U,I[p^{r+1}]) \rightarrow H^0(U,I[p^r])$ are surjective,
and as in the proof of Lemma 10 we conclude that the corestrictions $cor: H^0(U',I[p^n]) \rightarrow H^0(U,I[p^n])$
are surjective for open subgroups $U' \subset U$ of $G$. Therefore the system
$I[p^n]$ is acyclic for $\lim\limits_{\leftarrow, U,n}$, noting that
$$
\lim\limits_{\leftarrow, U, n} H^0(U,I[p^n]) = \lim\limits_{\leftarrow, n} H^0(U_n,I[p^n])
$$
for a cofinal family $(U_n)$ of subgroups between $H$ and $G$.

Therefore we get a spectral sequence
$$
E_2^{p,q}= H^p_{\mbox{\scriptsize cont}} (G,\, H;\, R^q\underline
T_pA)\Rightarrow H^{p+q}_{\mbox{\scriptsize cont}} (G,\, H;\,
R\underline T_pA)\; .
$$
From the snake lemma one immediately gets
$$
R^q\underline T_pA=\left\{\begin{array}{cc} (A/p^mA) & q=1\\
0 & q>1
\end{array}\right.
$$
(note that the described functor $A \rightsquigarrow (A/p^m)$ is effacable, since $A$ embeds into
an injective, hence divisible $G$-module. Hence we get a short exact sequences
$$
0\rightarrow H^n_{\mbox{\scriptsize cont}}(G,\, H;\, \underline
T_pA)\rightarrow H^n_{\mbox{\scriptsize cont}}(G,\, H;\,
R\underline T_pA)\rightarrow H^{n-1}_{\mbox{\scriptsize cont}}
(G,\, H;\, R^1\underline T_pA)\rightarrow 0 \, .
$$
This shows (a) and also explains the notation for $R^nF$. In fact,
$H^n_{\mbox{\scriptsize cont}}(G,\, H;\, R$ $\underline T_pA)$ is the
hypercohomology with respect to $H^0_{\mbox{\scriptsize
cont}}(G,\, H;\, -)$ of a complex $R\underline T_pA$ in
$M_G^{\mathbb N}$ computing the $R^i\underline T_pA$.

\ms\noi (b): If $H$ is a normal subgroup, we can regard the
functor $F$ as a functor from $M_G$ to the category $Mod_{\Lambda(
G/H)}$ of $\Lambda(G/H)$-modules. On the other hand, we can also
write $F$ as the composition of the left exact functors
$$
H^0(H,\, -):~M_G  \rightarrow  M_{G/H}~~~,~~~A  \rightsquigarrow
A^H
$$
and
$$
\tilde{F}:M_{G/H} \!\rightarrow  Mod_{\Lambda(G/H)},B
\shortmid\hspace{-0,8mm}\rightsquigarrow
\lim\limits_{\stackrel {\displaystyle\leftarrow} {n}}\lim\limits_{\stackrel {\displaystyle\leftarrow} {U/H}} H^0(U/H, B[p^n]) \!= H^0(G/H,\{1\};RT_pB) .
$$
(Note that $U/H$ runs through all open (normal) subgroups of
$G/H$.) This immediately gives the spectral sequence in (b).

\ms\noi (c): We claim that the functor $F$ is isomorphic to the
functor
$$
\begin{array}{rcl}
M_G & \rightarrow & Mod_{\Lambda(G/H)}\\
B & \shortmid\hspace{-0,8mm}\rightsquigarrow & \mbox{Hom
}_{\Lambda(G)}\,(B^\vee,\, \Lambda(G/H)).
\end{array}
$$
In fact, writing $\mbox{Hom}_{\Lambda(G)}(-,-)$ for the
homomorphism groups of compact $\Lambda(G)$-modules, we have (cf.
[Ja 3] p. 179)
$$
\begin{array}{rcl}
\mbox{Hom}_{\Lambda(G)}(B^\vee,\, \Lambda(G/H)) & = &
\lim\limits_{\stackrel {\displaystyle\leftarrow} {U}} \mbox{Hom}_{\Lambda(G)}(B^\vee,\, \mathbb Z_p[G/U])\\
& = &
\lim\limits_{\stackrel {\displaystyle\leftarrow} {n}}\lim\limits_{\stackrel {\displaystyle\leftarrow} {U}}\mbox{Hom}_{cont}(H^0(U,B)^\vee,\,\mathbb Z/p^n\mathbb
Z)\\
& = &
\lim\limits_{\stackrel {\displaystyle\leftarrow} {n}}\lim\limits_{\stackrel {\displaystyle\leftarrow} {U}}\mbox{Hom}_{cont}
(H^0(U,B[p^n])^\vee,\,\mathbb Z/p^n\mathbb Z)\\
& = &
\lim\limits_{\stackrel {\displaystyle\leftarrow} {n}}\lim\limits_{\stackrel {\displaystyle\leftarrow} {U}}H^0(U,B[p^n])\, ,
\end{array}
$$
where $U$ runs through all open subgroups of $G$ containing $H$,
and hence
$$
\mbox{Hom}_{\Lambda(G)}(B^\vee,\, \Lambda(G/H)) = H^0(G,H; RT_pB)
.
$$
Since taking Pontrjagin duals is an exact functor $M_G\rightarrow
C_G$ taking injectives to projectives, the derived functors of the
functor $B ~~\rightsquigarrow~~\mbox{Hom }_{\Lambda(G)}(B^\vee,\,
\Lambda(G/H))$ are the functors $B~~\rightsquigarrow~~
\mbox{Ext}^i_{\Lambda(G)}(B^\vee,\, \Lambda(G/H))$, and we get
(c). Finally, by applying (b) for $H' = H$ and (c) for $H = \{1\}$
we get (d).

\bs\noi Let us note that the proof of theorem 11 gives the
following $\mathbb Z/p^n$-analogue (by `omitting the inverse
limits over $n$'). For a profinite group $G$ let $\Lambda_n(G) =
\Lambda(G)/p^n = \mathbb Z/p^n[[G]]$ be the completed group ring
over $\mathbb Z/p^n$.

\bs\noi\textbf{Theorem 12} \emph{Let $H$ and $H'$ be normal
subgroups of a profinite group $G$, with $H' \subset H$, and let
$A$ be a discrete $\Lambda_n(G)$-module.\\
 \noi (a)
 There is a spectral sequence of $\Lambda_n(G/H)$-modules
 $$
 E_2^{p,q} = H^p(G/H', H/H';\, H^q(H',\, A)) ~\Rightarrow~
 H^{p+q} \,(G,\, H;\, A)\, .
 $$
 \noi (b) On the category of discrete $\Lambda_n(G)$-modules the
 $\delta$-functor $A \rightsquigarrow (H^m(G,$ $H; A)\mid m\geq 0)$ with values in
 the category of $\Lambda_n(G/H)$-modules is
 canonically isomorphic to the
 $\delta$-functor $A \rightsquigarrow (Ext^m_{\Lambda_n(G)}(A^\vee,\Lambda_n(G/H))$ $\mid m\geq
 0)$, where the $Ext$-groups are taken in the category of compact
 $\Lambda_n(G)$-modules.\\
 \noi (c) In particular, if $~\mathcal G = G/H$, and $A$ is a discrete
 $\Lambda_n(G)$-module, then one has a spectral sequence of
 $\Lambda_n(\mathcal G)$-modules}
 $$
 E_2^{p,q}\! = Ext^p_{\Lambda_n(\mathcal G)}(H^q(H,A)^\vee,
 \Lambda_n(\mathcal G))\Rightarrow H^{p+q} (G, H;
 A) = Ext^{p+q}_{\Lambda_n(G)}(A^\vee,\Lambda(\mathcal G)) .
 $$

\bs\noi\textbf{Corollary 13} \emph{With the notations as for
Theorem 1, let $A$ be a finite $\Lambda_n(G_S)$-module, and
$\Lambda_n = \Lambda(\mathcal G)$. Then there is a spectral
sequence of finitely generated $\Lambda_n$-modules
$$
E_2^{p,q} \!= \!Ext^p_{\Lambda_n}(H^q(G_{\infty,S},A)^\vee,
\Lambda_n)\!\Rightarrow\! \lim\limits_{\stackrel {\displaystyle\leftarrow} {k'}} H^{p+q}(G_S(k'), A)\! =\!
Ext^{p+q}_{\Lambda_n(G_S)}(A^\vee,\Lambda_n),
$$
where $k'$ runs through the finite subextensions $k'/k$ of
$k_\infty/k$.}

\bs\noi On the other hand, Theorem 1 also has the following
counterpart for finite modules.

\bs\noi\textbf{Theorem 14} \emph{With notations as for Theorem 1,
let $A$ be a finite $p$-primary $G_S$-module, of exponent $p^n$.
Then there is a spectral sequence
$$
E_2^{p,q}\! = \!Ext^p_\Lambda(H^q(G_{\infty,S},A)^\vee,
\Lambda)\!\Rightarrow \!\lim\limits_{\stackrel {\displaystyle\leftarrow} {k'}} H^{p+q-1}(G_S(k'), A)\! =
Ext^{p+q-1}_{\Lambda_n(G_S)}(A^\vee,\Lambda_n) ,
$$
where, in the inverse limit, $k'$ runs through the finite
extension $k'$ of $k$ inside $k_\infty$ and the transition maps
are the corestrictions.}

\bs\noi\textbf{Proof} As in the proof of Theorem 1, Theorem 11 (d)
applies to $G = G_S$ and $H = G_{\infty,S}$. But now the inverse
system $( A[p^n] )$ is Mittag-Leffler-zero in the sense of \cite{Ja1}:
if the exponent of $A$ is $p^d$, then the transition maps
$A[p^{n+d}] \rightarrow A[p^n]$ are zero. This implies that
$H^m_{cont}(G_S,(A[p^n])) = 0$ for all $m\geq 0$, cf. [Ja 1]. On
the other hand it is clear that the system $(A/p^n)$ is
essentially constant ($A/p^n = A$ for $n \geq d$). From Theorem 11
(a) and Lemma 10 we immediately get
$$
H^m_{cont}(G_S,G_{\infty,S};RT_pA) \cong
H^{m-1}_{cont}(G_S,G_{\infty,S};(A/p^n)) \cong
 \lim\limits_{\stackrel {\displaystyle\leftarrow} {k'}} H^{p+q-1}(G_S(k'),\, A)\, ,
$$
and hence the claim, by applying Theorem 12 (b) in addition.

\vspace{0.5cm}

\noi
Uwe Jannsen\\
Fakult\"at f\"ur Mathematik\\
Universit\"at Regensburg\\
93040 Regensburg\\
GERMANY\\
uwe.jannsen@mathematik.uni-regensburg.de\\
http://www.uni-regensburg.de/Fakultaeten/nat Fak I/Jannsen/index.html


\begin{thebibliography}{3}

\bibitem{Bi}
P. Billot, Quelques aspects de la descente sur une courbe elliptique dans le cas de r\'eduction supersinguli\`ere,
Compositio Math. {\bf 58} (1986), no. 3, 341--369.

\bibitem{Iw}
K. Iwasawa, On $\ZZ_\ell$-extensions of algebraic number fields. Ann. of Math. (2) 98 (1973), 246--326.

\bibitem{Ja1}
U. Jannsen, Continuous \'etale cohomology, Math. Ann. {\bf 280} (1988), 207--245.

\bibitem{Ja2}
U. Jannsen, On the $l$-adic cohomology of varieties over number fields and its Galois cohomology,
in {\it Galois groups over $ Q$ (Berkeley, CA, 1987)}, 315--360, Math. Sci. Res. Inst. Publ., 16,
Springer, New York, 1989

\bibitem{Ja3}
U. Jannsen, Iwasawa modules up to isomorphism, in {\it Algebraic Number theory}, 171--207, Adv. Stud. Pure Math. {\bf 17}, Academic
Press Boston 1989

\bibitem{La}
M. Lazard, Groupes analytiques $p$-adiques, Publ. Math. IHES {\bf 26} (1965), 389--603.

\bibitem{NSW}
J. Neukirch, A. Schmidt, K. Wingberg,
Cohomology of number fields.
Second edition. Grundlehren der Mathematischen Wissenschaften, 323. Springer-Verlag, Berlin, 2008. xvi+825 pp.

\bibitem{OV}
Y. Ochi, O. Venjakob,
On the structure of Selmer groups over p-adic Lie extensions.
J. Algebraic Geom. 11 (2002), no. 3, 547–580.

\bibitem{P-R}
B. Perrin-Riou, Arithm\'etique des courbes elliptiques et th\'eorie d'Iwasawa.
Mém. Soc. Math. France (N.S.) No. {\bf 17} (1984), 130 pp.

\bibitem{Sha}
S. Shatz, Profinite groups, arithmetic, and geometry.
Annals of Mathematics Studies, No. 67. Princeton University Press, Princeton, N.J.; University of Tokyo Press, Tokyo, 1972. x+252 pp.

\end{thebibliography}
\end{document}